\documentclass[11pt]{article}
\begin{document}
\bibliographystyle{plain}

\newtheorem{theorem}{Theorem}
\newtheorem{conjecture}{Conjecture}
\newtheorem{lemma}[theorem]{Lemma}
\newtheorem{note}[theorem]{Note}
\newtheorem{definition}{Definition}
\font\boldsets=msbm10
\def\F{{\hbox{\boldsets  F}}}
\def\N{{\hbox{\boldsets  N}}}
\def\R{{\hbox{\boldsets  R}}}
\def\Z{{\hbox{\boldsets  Z}}}
\def\blackslug{\hbox{\hskip 1pt \vrule width 4pt height 8pt depth 1.5pt
\hskip 1pt}}

\newcommand {\half}{\frac 1 2}
\newcommand{\ignore}[1]{}

\newcommand{\bbox}{\vrule height7pt width4pt depth1pt}
\def\QED{\quad\blackslug\lower 8.5pt\null\par}
\def\proof{\par\penalty-1000\vskip .5 pt\noindent{\bf Proof\/: }}
\def\prooft{\par\penalty-1000\vskip .5 pt\noindent{\bf Proof\/ }}
\newcommand{\prb}{\begin{problem}}
\newcommand{\eprb}{\bbox\end{problem}}
\newcommand {\upint}[1]   {\lceil {#1} \rceil}
\title
{Tales of Hoffman}
\author{
Yonatan Bilu
\thanks{
Institute of Computer Science,
Hebrew University Jerusalem 91904 Israel
\mbox{\em johnblue@cs.huji.ac.il}. 
This research is supported by the Israeli Ministry of Science. 
}
}
\date{}
\maketitle

\begin{abstract}
Hofmman's bound on the chromatic number of a graph states that 
$\chi \geq 1 - \frac {\lambda_1} {\lambda_n}$. Here we show that the
same bound, or slight modifications of it, hold for several graph 
parameters related to the chromatic number: the vector coloring number,
the $\psi$-covering number and the $\lambda$-clustering number.
\end{abstract}

\section{Introduction}

{\bf Hoffman's bound:}
Let $G$ be a graph on $n$ vertices, 
$\chi$ its chromatic number, and $A$ its adjacency matrix.
Let $\lambda_1$ 
and $\lambda_n$ the largest and least eigenvalues 
of $A$. A theorem of Hoffman \cite{Hoff} states that:
\[\chi \geq 1 - \frac {\lambda_1} {\lambda_n}.\]

{\bf The vector chromatic number:}
Karger, Motwani and Sudan \cite{KMS}
define a quadratic programming relaxation of the chromatic number,
called the {\em vector chromatic number}. This is 
the minimal $k$ such that there
exist unit vectors $u_1,\ldots,u_n \in \R^n$ with: 
\[<u_i,u_j> \leq - \frac 1 {k-1},\]
whenever $(i,j)$ is an edge in the graph.

Let $\chi_v$ denote the vector chromatic number of $G$. 
Karger, Motwani and Sudan observe that $\chi_v \leq \chi$.
In this note we show that Hoffman's bound holds for this parameter
as well:

\begin{theorem}\label{thrm-vec-color}
Let $G$ be a graph. Let $W \neq 0$ be a symmetric matrix such that
$W_{i,j} = 0$ whenever $(i,j) \not \in E$. 
Let $\lambda_1$ and $\lambda_n$ be the
largest and least eigenvalues of $W$. Then 
\[\chi_v(G) \geq 1 - \frac {\lambda_1} {\lambda_n}.\]
\end{theorem}

{\bf The $\psi$-covering number:}
Let $\psi$ be a graph parameter, such that 
$\psi = 1$ on graphs with no edges, 
and for every graph $G$, $\psi(G) \geq \chi(G)$.
The {\em $\psi$-covering} number of a graph $G$ was defined
by Amit, Linial and Matou\v{s}ek \cite{AmLiMa} 
to be the minimal $k$ such that there
exist $k$ subsets of $G$, $S_1,\ldots,S_k$ so that for every $v \in G$,
$\sum_{i\;:\;v \in S_i} \frac 1 {\psi(G[S_i])} \geq 1$.\\
They show that this value is bounded between 
$\sqrt{\chi(G)}$ and $\chi(G)$, and
ask whether better lower bounds can be proven
when $\psi(G) = dgn(G)+1$ (the degeneracy of
$G$ + 1), and $\psi(G) = \Delta(G)+1$ (the maximal 
degree in $G$ + 1).\\
To state our result, we'll need a couple of ad-hoc definitions:
\begin{definition}
A graph $G$ has {\em $c$-vertex cover} 
if there exists a cover $E(G) = \cup_{i\in V(G)} E_i$ such that for all 
$i \in V(G)$, $E_i \subset \{e \in E\;:\; i \in e\}$, and
$|E_i| \leq c$.\\
Denote by $L_{\psi,\alpha}(G)$ 
the minimal $k$ such that there
exist $k$ subsets of $G$, $S_1,\ldots,S_k$, so that 
for all $i=1,\ldots,k$, $G[S_i]$ has a $\alpha \cdot \psi(G)$-vertex cover, and
for every $v \in G$,
$\sum_{i\;:\;v \in S_i} \frac 1 {\psi(G[S_i])} \geq 1$.
\end{definition}
Observe that all graphs have 
a $\psi(G)$-vertex cover for $\psi(G) = dgn(G)+1$, 
and a $\half \psi(G)$-vertex cover when $\psi(G) = \Delta(G)+1$.
So in these cases, $L_{\psi,1}$ and $L_{\psi,\half}$, respectively, 
are exactly the
$\psi$-covering numbers.\\
Note also that $\alpha = 0$ means that the $S_i$ are independent sets. Thus,
since $\psi = 1$ on such sets, $L_{\psi,0} = \chi$.
\begin{theorem}\label{thrm-covering}
\[L_{\psi,\alpha}(G) \geq \frac { d - \lambda_n} {2\alpha - \lambda_n},\]
where $\lambda_n$ is the least eigenvalue of $G$, and $d$ the average degree.
\end{theorem}
Note that when the graph is regular and $\alpha = 0$,
this is the same as Hoffman's
bound.\\
For random 
$d$-regular graphs, $|\lambda_n| = O(\sqrt{d})$ and 
$\chi = \Theta(\frac d {\log d})$. So in this case (if $\alpha$ is
taken small) the bound  
is slightly better than $\sqrt{\chi}$ mentioned above.

{\bf The $\lambda$-clustering number:}
Finally, we are interested in a graph parameter that has to do with how well 
a graph can be partitioned into sparse clusters:
\begin{definition}
Let $W$ be a weighted adjacency matrix of a graph $G$. 
A partition $V = \dot \cup_{i=1}^k C_i$
is a {\em $\lambda$-clustering} of $G$ into $k$ clusters if
\[\max_{i \in [k]} \lambda_1(C_i) \leq \lambda,\]
where $\lambda_1(C_i)$ is the largest eigenvalue of the (weighted)
sub-graph spanned by the vertices 
in $C_i$.\\
The $\lambda$-clustering number of $G$ is the minimal $k$ such that there
exists a $\lambda$-clustering of $G$ into $k$ clusters.
\end{definition}
It is not hard to see that the $0$-clustering number is identical to the
chromatic number.\\
We show that Hoffman's bound can also be extended to this graph parameter:

\begin{theorem}\label{thrm-delta-clust}
Let $W$ be a weighted adjacency matrix. 
Let $\lambda_1$ 
and $\lambda_n$ the largest and least eigenvalues 
of $W$. The $\lambda$-clustering number of the graph is at least:
\[\frac {\lambda_1 - \lambda_n} {\lambda - \lambda_n}.\]
\end{theorem}

\section{Vectorial characterization of the least eigenvalue}
All three results mentioned in the previous section rely on the
following observation:
\begin{lemma}\label{obser}
Let $A$ be a real symmetric matrix and $\lambda_n$ its least eigenvalue.
\begin{eqnarray}\label{obs}
\lambda_n = &\min& \frac {\sum_{i,j=1}^n A_{i,j}<v_i, v_j>} 
{\sum_{i=1}^n||v_i||_2^2}.
\end{eqnarray}
where the minimum is taken over all $v_1,\ldots,v_n \in \R^n$.
\end{lemma}
\proof
By the Rayleigh-Ritz characterization, $\lambda_n$ equals
\begin{eqnarray*}
&\min& \sum_{i,j}A_{i,j} x_i x_j\\
&\mbox{s.t. }& x \in \R^n\\
&& ||x||_2 = 1.
\end{eqnarray*}
Denote by $PSD_n$ the cone of $n \times n$ positive semi-definite matrices.
For each unit vector $x \in \R^n$, let $X$ be the matrix $X_{i,j} = x_ix_j$. 
This is a positive semi-definite matrix of rank $1$ and trace $1$, 
and all such matrices
are obtained in this way. Hence, $\lambda_n$ equals:
\begin{eqnarray*}
&\min& \sum_{i,j} A_{i,j} X_{i,j}\\
&\mbox{s.t. }& X \in PSD_n\\
&& rank(X) = 1\\
&& tr(X) = 1.
\end{eqnarray*}
However, the rank restriction is superfluous. It restricts the solution to
an extreme ray of the cone $PSD_n$, but, by convexity, the optimum is attained
on an extreme ray anyway. Hence, $\lambda_n$ equals:
\begin{eqnarray*}
&\min& \sum_{i,j}A_{i,j}X_{i,j}\\
&\mbox{s.t. }& X \in PSD_n\\
&& tr(X) = 1.
\end{eqnarray*}
Now, think of each $X \in PSD_n$ as a Gram matrix of $n$ vectors,
$v_1,\ldots,v_n$ (i.e. $X_{i,j} = <v_i,v_j>$).
An equivalent formulation of the above is thus:
\begin{eqnarray*}
&\min& \sum_{i,j}A_{i,j}<v_i, v_j>\\
&\mbox{s.t. }& v_i \in \R^n \hspace{0.4in} \mbox{ for } i=1,\ldots,n\\
&& \sum_{i=1}^n||v_i||_2^2 = 1.
\end{eqnarray*}
Clearly, this is equivalent to $\ref{obs}$.
\QED

\section{Proofs of the theorems}
\prooft (Theorem \ref{thrm-vec-color}):
Let $G$ be a graph on $n$ vertices
with vector chromatic number
$\chi_v$. 
Let $W \neq 0$ be a symmetric matrix such that
$W_{i,j} = 0$ whenever $(i,j) \not \in E$. 
Let $\lambda_1$ and $\lambda_n$ be the
largest and least eigenvalues of $W$.\\
We choose vectors $v_1,\ldots,v_n$, and look at the bound
they give on $\lambda_n$ in Lemma \ref{obser}. 
Let $u_1, \ldots, u_n \in S^n$ be vectors on which the vector chromatic number
is attained.
That is, $<u_i,u_j> \leq - \frac 1 {\chi_v-1}$ for $(i,j) \in E$,
and $||u_i||_2 = 1$. Let $\alpha \in \R^n$ be an eigenvector of $W$
corresponding to $\lambda_1$.
Set $v_i = \alpha_i \cdot u_i$.\\
Since $W_{i,j} = 0$ whenever $<u_i,u_j> \;> - \frac 1 {\chi_v-1}$, by 
Lemma \ref{obser}, 
\begin{eqnarray*}
\lambda_n &\leq& 
\frac {\sum_{i,j}W_{i,j} \alpha_i \alpha_j <u_i, u_j>} 
{\sum_{i=1}^n\alpha_i^2 \cdot ||u_i||_2^2} \leq
- \frac 1 {\chi_v-1} \cdot \frac {\sum_{i,j} W_{i,j} \alpha_i \alpha_j}
{\sum_i \alpha_i^2} = \\
&=&
- \frac 1 {\chi_v-1} \cdot
\frac {\alpha^t W \alpha} {||\alpha||^2}
= - \frac 1 {\chi_v-1} \cdot \lambda_1.
\end{eqnarray*} 
Equivalently, $\chi_v \geq 1 - \frac {\lambda_1} {\lambda_n}$, as claimed.
\QED

\prooft (Theorem \ref{thrm-covering}):
Denote $k = L_{\psi,\alpha}(G)$, and let 
$u_1,\ldots,u_k$ be the vertices
of the regular $(k-1)$-dimensional 
simplex centered at $0$ - i.e. $<u_i,u_j> = 1$
when $i=j$ and $ \frac {-1} {k-1}$ otherwise. 
Again we choose vectors $v_1,\ldots,v_n$. We do so probabilistically.
Let $S_1,\ldots,S_k$ be the subsets attaining the value $k$. 
For each $i$, $v_i$ will be chosen from among
the $u_j$'s such that $i \in S_j$. 
Specifically, let $h_i = \sum_{j\;:\;i \in S_j}\frac 1 {\psi(S_j)}$.
The probability that $v_i$ is chosen to be $u_j$ is 
$p_{i,j} = h_i^{-1}\frac 1 {\psi(S_j)}$.  
Note that $h_i \geq 1$, and so $p_{i,j} \leq \frac 1 {\psi(S_j)}$.
(there is a slight abuse of notation here - by $\psi(S_j)$ we refer to
$\psi(G[S_j])$.)\\
Say that an edge is ``bad'' if both its endpoints are assigned the same
vector. For a given $j$, the probability that an edge $(i,i') \in E(S_j)$
is ``bad'' because
both endpoints were assigned to $u_j$ is $p_{i,j} p_{i',j}$. Thus, the 
expected number of ``bad'' edges is at most:
\begin{eqnarray*}
\sum_{j=1}^k \sum_{(i,i')\in E(S_j)}p_{i,j} p_{i',j}.  
\end{eqnarray*}
Each $S_j$ has a $\alpha \cdot \psi(G)$-vertex cover 
$E(S_j) = \cup E_j^i$. Summing 
the expression above according to this cover (some edges might be counted
more than once) we get that the expected number of ``bad'' edges is at most:
\begin{eqnarray*}
\sum_{j=1}^k \sum_{i \in S_j} 
\sum_{i'\;:\;(i,i')\in E_j^i} p_{i,j} p_{i',j} \;&\leq&\;
\sum_{j=1}^k \sum_{i \in S_j} 
p_{i,j} |E_j^i| \frac 1 {\psi(S_j)} \;\leq\;
\sum_{j=1}^k \sum_{i \in S_j} 
p_{i,j} \alpha\\ &=& 
\alpha \sum_{i \in G} \sum_{j\;:\;i \in S_j} p_{i,j} = \alpha n.
\end{eqnarray*}
In particular, there is a choice of $v_i$'s such that the number of ``bad''
edges is at most $\alpha n$. Assume this is the case.
If $(i,j)\in E(G)$ is a ``bad'' edge then ${<v_i,v_j> = 1}$. Otherwise
\mbox{$<v_i,v_j> = - \frac 1 {k-1}$}.\\
Lemma \ref{obser} now gives:
\begin{eqnarray*}
\lambda_n \leq \left(2\alpha n - \frac 1 {k-1} (dn - 2\alpha n)\right) / n =
\frac {2k\alpha} {k-1} - \frac d {k-1},
\end{eqnarray*}
or $(k-1) \lambda_n \leq 2k \alpha - d$. Equivalently, 
$k \geq \frac { d - \lambda_n} {2\alpha - \lambda_n}$.
\QED
\begin{note}
In the definition of the $\psi$-covering number, and of
$L_{\psi,\alpha}$, it is required 
that for every $v \in G$, 
$\sum_{i\;:\;v \in S_i} \frac 1 {\psi(S_i)} \geq 1$. 
The bounds given in \cite{AmLiMa} hold also 
if we demand that all sums {\em equal}
$1$. In this case, the theorem holds also if we relax the condition that
all $S_i$ have an $\alpha \cdot \psi(G)$-vertex cover, 
and require only that for each $S_i$,
$|E(S_i)| \leq \alpha \cdot \psi(S_i) \cdot |V(S_i)|$.
\end{note}

\prooft (Theorem \ref{thrm-delta-clust}):
Denote the $\lambda$-clustering number of $W$ by $k$.
Let $u_1, \ldots, u_k \in R^n$ be the vertices of a regular simplex centered 
at the origin, as above.
Let $\alpha \in \R^n$ be an eigenvector of $W$, 
corresponding to $\lambda_1$. Let $C_1,\ldots,C_k$ be a $\lambda$-clustering
of $G$. Define $\phi:V \rightarrow [k]$ to be the index of the cluster
containing a vertex. That is, $i \in C_{\phi(i)}$. 
Define $W_t$ to be the weighted
adjacency matrix of the sub-graph spanned by $C_t$ 
(So $\lambda_1(W_t) \leq \lambda$).
Set $v_i = \alpha_i \cdot u_{\phi(i)}$.\\
By Lemma \ref{obser}, 
\begin{eqnarray*}
\lambda_n &\leq& 
\frac {\sum_{i,j}\alpha_i \alpha_j <u_{\phi(j)}, u_{\phi(j)}> W_{i,j}} 
{\sum_{i=1}^n\alpha_i^2 \cdot ||u_{\phi(i)}||_2^2}\\
&=&
- \frac 1 {k-1} \cdot \frac 
{\sum_{i,j \;:\; \phi(i) \neq \phi(j)}\alpha_i \alpha_j W_{i,j}}
{\sum_i \alpha_i^2} + \frac
{\sum_{i,j \;:\; \phi(i) = \phi(j)}\alpha_i \alpha_j W_{i,j}}
{\sum_i \alpha_i^2} 
\\
&=&
- \frac 1 {k-1} \cdot \frac {\alpha^t W \alpha} {||\alpha||^2} +
\frac k {k-1}
\frac {\sum_{t=1}^k\sum_{i,j \in C_t}\alpha^t W_t \alpha}
{\sum_i \alpha_i^2} \\
&\leq&
- \frac 1 {k-1} \lambda_1 + \frac k {k-1} \lambda
\end{eqnarray*} 
Equivalently, $k \geq \frac {\lambda_n - \lambda_1} {\lambda_n - \lambda}$
\QED

\section*{Acknowledgments}
I thank Nati Linial, Eyal Rozenman and Yael Vinner
for their contribution to this note.

\bibliography{bib}

\begin{thebibliography}{1}

\bibitem{AmLiMa}
A.~Amit, N.~Linial, and J.~Matou{\v{s}}ek.
\newblock Random lifts of graphs: independence and chromatic number.
\newblock {\em Random Structures Algorithms}, 20(1):1--22, 2002.

\bibitem{Hoff}
A.~J. Hoffman.
\newblock On eigenvalues and colorings of graphs.
\newblock In {\em Graph Theory and its Applications (Proc. Advanced Sem., Math.
  Research Center, Univ. of Wisconsin, Madison, Wis., 1969)}, pages 79--91.
  Academic Press, New York, 1970.

\bibitem{KMS}
D.~Karger, R.~Motwani, and M.~Sudan.
\newblock Approximate graph coloring by semidefinite programming.
\newblock {\em J. ACM}, 45(2):246--265, 1998.

\end{thebibliography}
\end{document}